\documentclass[16pt]{article}
\usepackage[pdftex,bookmarksopen=true,bookmarks=true,unicode,setpagesize]{hyperref}
\hypersetup{colorlinks=true,linkcolor=black,citecolor=black}

\usepackage{amssymb, amsmath, amsthm}

\usepackage{nicefrac}
\usepackage[symbol]{footmisc}
\usepackage{graphicx}
\usepackage{setspace}
\usepackage{tocloft}
\usepackage{lineno}

\usepackage{xcolor}

\allowdisplaybreaks
\usepackage{cite}

\newcommand{\vertiii}[1]{{\vert\kern-0.25ex\vert\kern-0.25ex\vert #1
		\vert\kern-0.25ex\vert\kern-0.25ex\vert}}

\theoremstyle{remark}

\theoremstyle{remark}

\theoremstyle{remark}

\begin{document}

	\vspace{-20mm}
	{\Large \bf{Diagonal Condition in Multiplication Table of $\, \mathbb{Z} [i] / (\alpha)$}}

	\begin{center}
	{\large Chadaphorn Kodsueb}\\ \vspace{1mm} School of Mathematical Sciences and Geoinformatics, Institute of Science, \\ Suranaree University of Technology, Nakhon Ratchasima, Thailand 30000; \\ e-mail: \texttt{Chadaphorn.Ko@sut.ac.th}\vspace{2mm}
	\end{center}

	{\small
		\begin{center}
			{\bf Abstract}
		\end{center}
			\qquad Multiplication table of a ring with identity $\mathbf{1}$ is said to have the diagonal condition if $\mathbf{1}$’s occur only on the main diagonal. In this paper, we study the diagonal condition in the ring of Gaussian integers $ \mathbb{Z} [i]$. Furthermore, we also find the Gaussian integers $\alpha$ so that the rings of Gaussian integers modulo $\alpha$ satisfy the diagonal condition.
		\noindent   
		\noindent
		
	} \vspace{2mm}
	
	{\bf Keywords:}  diagonal condition, multiplication table, Gaussian integers, quotient ring
	\vspace{2mm}
	
	{\bf 2020 MSC:} 	11R04, 11T30, 11N13.
	
	\section{Introduction}
	
	\qquad In \cite{Chebolu-24} S.K. Chebolu (2012), he observed that $\mathbf{1}$'s only appear on the main diagonal of multiplication tables of $\mathbb{Z}_2$, $\mathbb{Z}_3$, $\mathbb{Z}_4$ and $\mathbb{Z}_8$ but not for $\mathbb{Z}_5$. The aim of \cite{Chebolu-24} S.K. Chebolu (2012) is to find which positive integers $n$ that make $\mathbf{1}$'s occur only on the main diagonal (not off diagonal) of multiplication table of $\mathbb{Z}_n$. Given $\mathbb{Z}_n = \{ 1,2, \ldots, n-1 \}$. Any ring which has $\mathbf{1}$'s only appear on the main diagonal of its multiplication table. One says the ring satisfies \textit{diagonal condition}, summarized in the following. 
	
	\vspace{2ex}\noindent{\textbf{Theorem 1.1} \ Let $n$ be a positive integer for $\mathbb{Z}_n$. The statements are equivalent:}
	\begin{enumerate}   
		\item $\mathbf{1}$'s occur only on the main diagonal (not off diagonal) of multiplication table of $\mathbb{Z}_n$.
		\item If $a$ is an element which has multiplicative inverse in $\mathbb{Z}_n$, then $a^2 = 1 \in \mathbb{Z}_n$.
		\item If $a$ is an integer relatively prime with $n$, then $n \mid a^2 - 1$.
		\item If $p$ is a prime such that $p \nmid n$, then $n \mid p^2 - 1$.
	\end{enumerate}
	
	Finally, he found that $\mathbf{1}$’s occur only on the main diagonal of multiplication table of $\mathbb{Z}_n$ if and only if $n$ is a divisor of $24$. This became the main theorem in \cite{Chebolu-24} S.K. Chebolu (2012). He also provided five distinct proofs for this theorem. However, the special property about divisors of $24$ has not ended there. In 2024, J. C. Lagarias proposed problems on Floor Quotients in American Math. Monthly problem 11747, Vol. 121, p. 83. His solution appears in Vol. 123 (2016), pp. 199--200\cite{fromJeff-2}.
	
	Later, \cite{Chebolu-12} Chebolu, S.K. and Mayers (2013) investigated on the (same) diagonal condition of multiplication table but in a different ring, $\mathbb{Z}_n [x]$. They had very interest result, and then formed it as a main theorem ``For any positive integer $m$, the multiplication table for the polynomial ring $\, \mathbb{Z}_n [x_1, x_2, \cdots, x_m]$ has $\mathbf{1}$’s only on the diagonal if and only if n is a divisor of 12". Another result involved with number $12$ has shown in \cite{fromJeff-3} Poonen, B. and Rodriguez-Villegas, F. (2000). They discussed the number $12$ related to convex lattice polygons with one interior lattice point.
	
	Additionally, there is one paper linked to both numbers $12$ and $24$. In \cite{fromJeff-1} Godinho, L., von Heymann, F. and Sabatini, S. (2017), they were interested in connection with mirror-symmetry and Calabi-Yau toric varieties. This kind of thing leads into aspects of Topology and Geometry.
	
	Beyond than that \cite{Sophie-Germain} Genzlinger, K. and Lockridge, K. (2015) studied the diagonal condition on multiplication cube of $\, \mathbb{Z}_n$. They gained the main theorem as ``All $\mathbf{1}$'s in the multiplication cube  for $\mathbb{Z}_n$ lie exclusively on the diagonal
	or the coordinate planes (where one of the three coordinates is $1$) if and only if  $n$ is a divisor of $4$ or $6$.
	
	Nevertheless, the author's aim is still about the diagonal condition in multiplication table of (another) classical ring: Gaussian integers $\mathbb Z [i]$, specifically a ring of Gaussian integers modulo $\alpha$ .

	\section{Gaussian integers} \label{GI}
	
	\qquad This section provides crucial background of Gaussian integers $\mathbb{Z}[i]$ and Gaussian integers modulo $\alpha$ where $\alpha$ is nonzero in $\mathbb{Z}[i]$. Some materials in this section, we take it partially from \cite{K.Conrad} K. Conrad (2016), \cite{Dummit-Foote} D.S. Dummit and R.M Foote (2004) and  \cite{ContempAA} J.A. Gallian (2021). 
	
	\subsection{Basic knowledge of $\mathbb{Z}[i]$}
	\label{Basic knowledge}
	
	\qquad Gaussian integers are complex numbers is in the form $a+bi$ where $a, b \in \mathbb{Z}$, denoted $\mathbb{Z} [i]$. In fact, $\mathbb{Z} [i]$ is subring of $\mathbb{C}$ with $1$ is multiplicative identity. Since $\mathbb{C}$ is a field, $\mathbb{Z} [i]$ is an integral domain.
	
	\vspace{2ex}\noindent{\textbf{Definition} \  For $\alpha = a+bi$, \textit{Norm} of $\alpha$ denote by $\mathcal{N} (\alpha) := \alpha \overline{\alpha} = (a+bi)(a-bi) = \alpha^2 + \beta^2 . $}
	
	\vspace{2ex}\noindent{\textit{Remark 2.1.1} \ It is obvious that $\, \mathcal{N} (\alpha) \in {\mathbb{N} }_0 = \mathbb{Z}^+ \cup \{ 0 \}$ if including the case $\alpha = 0$.} 
	
	\vspace{2ex}\noindent{\textbf{Theorem 2.1.2} \ For $\alpha, \beta \in \mathbb{Z} [i]$, \ $\mathcal{N} (\alpha\beta) = \mathcal{N} (\alpha) \, \mathcal{N} (\beta)$.}  
	
	\vspace{2ex}\noindent{\textbf{Definition} \ The Gaussian integer $\beta$ is a \textit{unit} if $\beta$ has multiplicative inverse in $\mathbb{Z} [i]$ i.e. there exists $\gamma \in \mathbb{Z} [i]$ such that $\gamma\beta = 1$.}
	
	\vspace{2ex}\noindent{\textbf{Theorem 2.1.3} \ The only Gaussian integers that are units in $\mathbb{Z} [i]$ are $1, -1, i, -i$. }
	
	\vspace{2ex}\noindent{\textbf{Definition} \  For $\alpha, \beta \in \mathbb{Z} [i]$ where $\beta \not= 0$, we say that $\beta \mid \alpha$ if $\alpha = \beta\gamma$ for some $\gamma \in \mathbb{Z} [i]$. Then $\beta$ is called a \textit{divisor} of $\alpha$. }
	
	\vspace{2ex}\noindent{\textbf{Theorem 2.1.4} \ For all $\alpha, \beta$ in $\mathbb{Z} [i]$, if $\beta \mid \alpha$ in $\mathbb{Z} [i]$ then $\mathcal{N} (\beta) \mid \mathcal{N} (\alpha)$ in $\mathbb{Z}$. }
	
	\vspace{2ex}\noindent{\textbf{Theorem 2.1.5} \big(Division Algorithm in $\mathbb{Z} [i]$\big) \ Given $\alpha, \beta \in \mathbb{Z} [i]$ with $\beta \not= 0$. There exist $\kappa, \rho $ in $\mathbb{Z} [i]$ such that $\alpha = \beta\kappa + \rho$ and $\mathcal{N} (\rho) < \mathcal{N} (\beta)$. In fact, we can choose $\rho$ that makes $\mathcal{N} (\rho) < \nicefrac{1}{2} \, \mathcal{N} (\beta)$. }
	\begin{proof} \ We provide proof without words in Geometry: Division Algorithm in $\mathbb{Z} [i]$ .
		
		\begin{center}
			\includegraphics[scale=0.80]{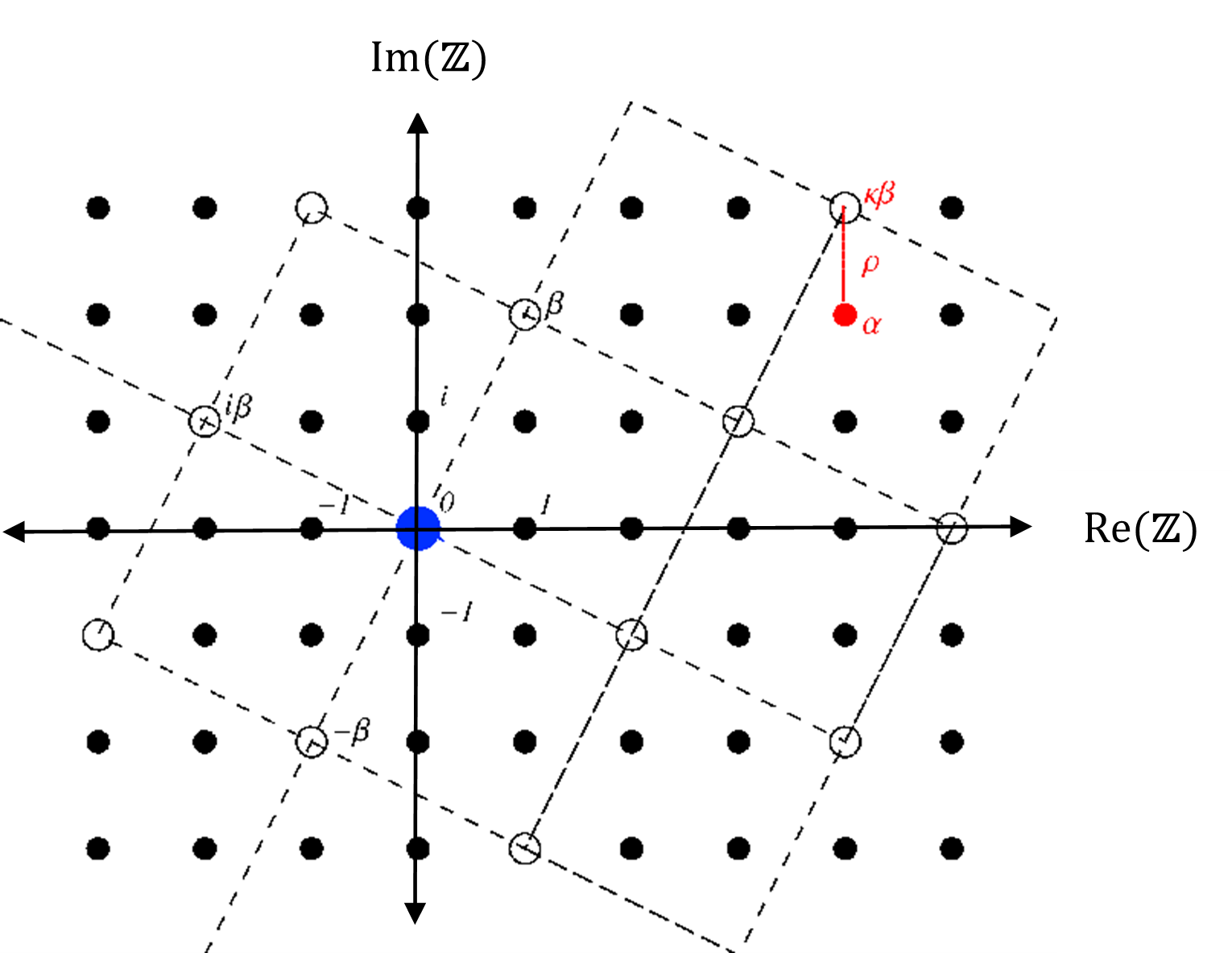}
		\end{center}
		
		\footnote[4]{Author's explanation}From the figure, vector from $0$ to $i\beta$ is perpendicular with vector from $0$ to $\beta$ in Complex Plane $\mathbb{C} = \mathbb{R}^2$. So the set $\beta \mathbb{Z} [i]  = \{ \beta\kappa \mid \beta \in \mathbb{Z} [i] \}$ forms squares with length $\lvert \beta \rvert = \sqrt{\mathcal{N}(\beta)}$ netting as lattice. It is easy to see that $\alpha$ appearing inside at least 1 square. Let $\beta\kappa$ be the most nearest corner point with point $\alpha$. Given $\rho = \alpha - \beta\kappa \in \mathbb{Z} [i]$. Then $\lvert \rho \rvert$ is less than or equal to half a diagonal length, i.e. $\sqrt{\mathcal{N}(\rho)} = \lvert \rho \rvert \leq \frac{\sqrt{2}}{2} \lvert \beta \rvert < \sqrt{\mathcal{N}(\beta)}$.
	\end{proof}
	
	Known fact is $\mathbb{Z} [i]$ is a Euclidean Domain (\textbf{ED}). Then $\mathbb{Z} [i]$ is also a \textbf{PID} and a \textbf{UFD}.
	
	\vspace{2ex}\noindent{\textbf{Definition} \ Let $\gamma \in \mathbb{Z} [i]$, we call $1, -1, i, -i, \gamma, -\gamma, \gamma i, -\gamma i$ are \textit{trivial factors} of $\gamma$.}
	
	\vspace{2ex}\noindent{\textbf{Definition} \ For $\alpha \in \mathbb{Z} [i] $ which is not a unit, $\alpha$ is a \textit{factorized element (\textit{reducible element})}. If $\alpha$ is not a trivial factor and $\alpha$ has only trivial factors, then $\alpha$ is a \textit{irreducible element} or called a \textit{Gaussian prime}.}
	
	\vspace{2ex}\noindent{\textbf{Theorem 2.1.6} \ If $\mathcal{N} (\alpha)$ is a prime in $\mathbb{Z}$, then $\alpha$ is a Gaussian prime in $\mathbb{Z} [i]$.  }
	
	\vspace{2ex}\noindent{\textbf{Definition} \ A \textit{common divisor} of $\alpha$ and $\beta$ is a Gaussian integer $\delta$ such that $\delta \mid \alpha$ and $\delta \mid \beta$. Denote $\delta$ := $\gcd(\alpha, \beta)$ a \textit{greatest common divisor} of $\alpha$ and $\beta$ if and only if 
		\begin{enumerate}
			\item $\delta \mid \alpha$ and $\delta \mid \beta$,
			\item For all $\gamma \in \mathbb{Z} [i]$, if $\gamma \mid \alpha$ and $\gamma \mid \beta$ then $\gamma \mid \delta$.    
	\end{enumerate} }
	
	One says $\alpha$ and $\beta$ are \textit{relatively (Gaussian) primes} if  $\delta$ = $\gcd(\alpha, \beta)$ is a unit in $\mathbb{Z} [i]$.
	
	\vspace{2ex}\noindent{\textbf{Theorem 2.1.7} \ Given $\pi$ is a Gaussian prime. For all $\alpha, \beta \in \mathbb{Z} [i]$, if $ \pi \mid \alpha\beta $ then $\pi \mid \alpha$ or $\pi \mid \beta$. }
	
	\vspace{2ex}\noindent{\textit{Remark 2.1.8} \ Since $\mathbb{Z} [i]$ is a Unique Factorization Domain, Gaussian integers which are nonzero and non-unit can be factorized as multiplication of distinct Gaussian primes and a unit. }
	
	\subsection{Coset representatives of quotient ring $\, \mathbb{Z}[i] \slash (\alpha)$}
	\label{subsect:Coset representatives}
	
	\qquad This subsection gives concepts of congruence--modulo and coset representatives in $\mathbb{Z} [i] \slash (\alpha)$.
	
	\vspace{2ex}\noindent{\textbf{Definition} \  For Gaussian integers $\alpha, \beta$ and $\gamma$, we say that $\beta$ \textit{congruence} $\gamma$ \textit{modulo} $\alpha$ denoted by $\beta \equiv \gamma \pmod{\alpha}$, when $\alpha \mid (\beta - \gamma)$. }
	
	\vspace{2ex}\noindent{\textit{Remark 2.2.1} \ It can be proved that Congruence--Modulo is an equivalence relation on $\, \mathbb{Z}[i]$. The equivalence classes are cosets in quotient ring  $\, \mathbb{Z}[i] \slash (\alpha)$ .}
	
	\vspace{2mm}
	
	Next, we provide methods to find coset representatives of $\, \mathbb{Z}[i] \slash (\alpha)$ when $\alpha = a+bi$ through multiples of Gaussian in grid as followed.
	\begin{enumerate}
		\item Let $m, n \in \mathbb{Z}$. Consider $(a+bi)(m+ni) = (a+bi)m + (-b+ai)n$.
		\item Plot points $(a,b)$ and $(-b,a)$ in the grid, then form two vectors from origin to these points.
		\item Extend the sided of two vectors to create a (main) square. Coset representatives of $\, \mathbb{Z}[i] \slash (\alpha)$ are points located inside the square and on two main vectors including origin.
		\item Translate the square till it can be filled the grid. Then, we have (all) cosets.
	\end{enumerate}
	
	\vspace{2ex}\noindent{\textit{Examples 2.2.2} \ Find coset representatives of $\, \mathbb{Z}[i] \slash (2+2i)$. } \vspace{1mm}
	
	First, consider $(2+2i)(m+ni) = (2+2i)m + (-2+2i)n$. Then follow steps 2. and 3. :  \vspace{2mm}
	\begin{center}
		\includegraphics[scale=0.80]{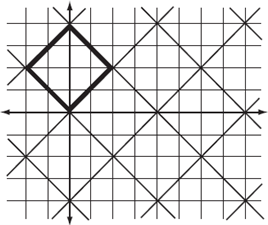}
	\end{center}
	\qquad From the figure, the coset representatives of $\, \mathbb{Z}[i] \slash (2+2i)$ are $0, i, 2i, 3i, 1+i, 1+2i, \\ -1+i, -1+2i$ equivalent with points $(0,0), (0,1), (0,2), (0,3), (1,1), (1,2), (-1,1), (-1,2)$ in grid.
	
	\vspace{2ex}\noindent{\textbf{Lemma 2.2.3} \ For $a, b, c, d, m \in \mathbb{Z}$, we have }
	\begin{enumerate}   
		\item $a \equiv b \pmod{c}$ in $\mathbb Z$ if and only if $a \equiv b \pmod{c}$ in $\mathbb{Z}[i]$.
		\item $a+bi \equiv c+di \pmod{m}$ if and only if $a \equiv c \pmod{m}$ and $b \equiv d \pmod{m}$.
	\end{enumerate}
	
	Let $\mathfrak{N} (\alpha)$ be the number of coset representatives of $\, \mathbb{Z}[i] \slash (\alpha)$ i.e. $\mathfrak{N} (\alpha) := \lvert \mathbb{Z}[i] \slash (\alpha) \rvert$, $\alpha \in \mathbb{Z}[i]$.
	
	\vspace{2ex}\noindent{\textbf{Lemma 2.2.4} \  If $m \in \mathbb{Z}$ is nonzero, then $\, \mathfrak{N} (m) = m^2$ .}
	
	\vspace{2ex}\noindent{\textbf{Lemma 2.2.5} \  If $\alpha \not= 0$ in $\mathbb{Z}[i]$, then $\, \mathfrak{N} (\overline{\alpha}) = \mathfrak{N} (\alpha)$ .}
	
	\vspace{2ex}\noindent{\textbf{Lemma 2.2.6} \  If $\alpha$ and $\beta$ are nonzero in $\mathbb{Z}[i]$, then $\, \mathfrak{N} (\alpha\beta) = \mathfrak{N} (\alpha) \, \mathfrak{N} (\beta)$. In another word, function $\mathfrak{N}$ is multiplicative.}
	
	\vspace{2ex}\noindent{\textbf{Theorem 2.2.7} \ If  $\alpha \not= 0$ in $\mathbb{Z}[i]$, then $\, \mathfrak{N} (\alpha) = \mathcal{N} (\alpha)$.}
	\begin{proof} \ 
		By Lemma~2.2.6, $\, \mathfrak{N} (\alpha \overline{\alpha}) = \mathfrak{N} (\alpha) \, \mathfrak{N} (\overline{\alpha})$. From Lemma~2.2.5, we have the right hand side of previous equation is $\, \mathfrak{N} (\alpha)^2$. Since $\alpha \overline{\alpha} = \mathfrak{N} (\alpha)$ is integer and Lemma~2.2.4, it can be concluded that $\, \mathfrak{N} (\alpha \overline{\alpha}) = \mathcal{N} (\alpha)^2$. Thus, $\, \mathcal{N} (\alpha)^2 = \mathfrak{N} (\alpha)^2$. Hence, we have $\, \mathcal{N} (\alpha) = \mathfrak{N} (\alpha)$.
	\end{proof}
	
	\subsection{Gaussian primes}
	\label{subsect:Gaussian primes}
	
	\qquad In this subsection, we find all Gaussian primes from Gaussian integers which divide primes in $\mathbb{Z}$.
	
	\vspace{2ex}\noindent{\textbf{Definition} \  In $\mathbb{Z}[i]$, we say that $\alpha \sim \beta$ if $\alpha = u\beta$ or $\alpha = u\beta$ where $u$ is a unit in $\mathbb{Z}[i]$. }
	
	\vspace{2ex}\noindent{\textbf{Lemma 2.3.1} \ Let $\pi$ be a Gaussian prime. Then, there exists a prime $p \in \mathbb{Z}$ such that $\, \pi \mid p$.} \vspace{2mm}
	
	We can find all Gaussian primes from Gaussian integers which divide primes in $\mathbb{Z}[i]$.
	
	\vspace{2ex}\noindent{\textbf{Theorem 2.3.2} \ A prime in $\mathbb{Z}$ is a factorized element (reducible element) in $\mathbb{Z}[i]$ if and only if $p$ can be written as a sum of two squares. }
	
	\vspace{2ex}\noindent{\textit{Remark 2.3.3} \  Any primes $p \in \mathbb{Z}$ which cannot be written as a sum of two squares, then $p$ is not a factorized element (reducible element). Thus, $p$ is a Gaussian primes. }
	
	\vspace{2ex}\noindent{\textbf{Corollary 2.3.4} \ If $p$ is a prime in $\mathbb{Z}$ such that $p \equiv 3 \pmod 4$, then $p$ is a Gaussian prime. }
	
	\vspace{2ex}\noindent{\textbf{Lemma 2.3.5} \  Given $p \in \mathbb{Z}$ is a prime. Either $p=2$ or $p \equiv 1 \pmod 4$ if and only if there exist integers $a, b$ such that $p = a^2 + b^2$. }
	
	\vspace{2ex}\noindent{\textbf{Theorem 2.3.6} \   Let $p$ be a prime in $\mathbb{Z}$  such that $p \equiv 1 \pmod 4$. Then $p = \pi \bar{\pi}$ where $\pi = a+bi$ and $\bar{\pi} = a-bi$ whilst $\pi \not\sim \bar{\pi}$. }
	\begin{proof} \
		By Lemma 2.3.4, there exist integers $a, b$ such that $p = a^2 + b^2$. So, $p = \pi \bar{\pi}$ where $\pi = a+bi$ and $\bar{\pi} = a-bi$. Assume that $a+bi \sim a-bi$. Then, there are 4 cases as followed \begin{itemize}
			\item \textit{case} $a+bi = a-bi$: we have $b=0$. Then $p = a^2$, a contradiction.
			\item \textit{case} $a+bi = -(a-bi)$: we have $a=0$. Then $p = b^2$, a contradiction.
			\item \textit{case} $a+bi = i(a-bi)$: we get $a = b$. Then $p = 2 a^2 = 2 b^2$, a contradiction.
			\item \textit{case} $a+bi = -i(a-bi)$: we get $a = -b$. Then $p = 2 b^2 = 2 a^2$, a contradiction.
		\end{itemize}
		
		Thus, $\pi = a+bi \not\sim a-bi = \bar{\pi}$.
	\end{proof}
	
	\vspace{2ex}\noindent{\textbf{Theorem 2.3.7}\footnote[5]{Main theorem for this section.} \  These are all Gaussian primes:  \begin{enumerate}
			\item 1+i,
			\item $p$ where $p$ are primes in $\mathbb Z$ s.t. $p \equiv 3 \pmod 4$,
			\item $\pi = a+bi$ and $ \bar{\pi} = a-bi$ when $N(\pi) = N(\bar{\pi}) = a^2 + b^2 = p$, and $p$ are primes in $\mathbb Z$ s.t. $p \equiv 1 \pmod 4$.
	\end{enumerate}}
	\begin{proof} \
		We can find all Gaussian primes which divide prime $p \in \mathbb{Z}$ from Lemma 2.3.1: \begin{enumerate}
			\item When $p=2$, we have $2 = (1+i)(1-i) = -i(1+i)^2$. Then $1+i$ is satisfied the condition i.e. $1+i \mid 2$ and $1+i \not\sim 1-i$,
			\item When $p \equiv 3 \pmod 4$, $p$ is a Gaussian prime by Corollary 2.3.4,
			\item When $p \equiv 1 \pmod 4$, by Theorem 2.3.6, $\pi$ and $\bar{\pi}$ are Gaussian primes where $p = \pi \bar{\pi}$ and $\pi \not\sim \bar{\pi}$.
		\end{enumerate}
		
		Thus, we have proved.
	\end{proof}
	
	
	\section{Main results}
	
	\label{sect:sections}
	\qquad From subsection~\ref{subsect:Gaussian primes} (Theorem~2.3.7), we know that all Gaussian primes are in the form
	\begin{enumerate}
		\item $1+i$,
		\item $p$ where $p$ are primes in $\mathbb Z$ s.t. $p \equiv 3 \pmod 4$,
		\item $\pi = a+bi$ and $\bar{\pi} = a-bi$ when $N(\pi) = N(\bar{\pi}) = a^2 + b^2 = p$, and $p$ are primes in $\mathbb Z$ s.t. $p \equiv 1 \pmod 4$.
	\end{enumerate}
	
	We also know that any Gaussian integer $\alpha \ne 0$ and $\alpha$ not a unit can be written as Gaussian primes factorization in the following theorem. \vspace{2mm}
	
	\vspace{0.5ex}\noindent{\textbf{Theorem 3.1} If $\alpha = u \, \pi_1^{c_1} \pi_2^{c_2} \cdots \pi_r^{c_r}$ where  $u$  is unit and $\pi_i$ are Gaussian primes such that $u \nsim \pi_i$, then $$\mathbb Z [i] / (\alpha) \cong \mathbb Z [i] / (u) \times \mathbb Z [i] / (\pi_1^{c_1}) \times \mathbb Z [i] / (\pi_2^{c_2}) \times \cdots \times \mathbb Z [i] / (\pi_r^{c_r}). $$}
	\begin{proof} \ 
		This is trivial from Chinese Remainder Theorem, see statement in \cite{BookofThomas} T. Koshy (2002).
	\end{proof} 
	
	\vspace{2ex}\noindent{\textbf{Theorem 3.2} \ $\mathbb Z [i] / ( (1+i)^n )$ do not satisfy diagonal condition when $n \geq 3$.}
	\begin{proof} \ 
		Suppose that $\mathbb Z [i] / ( (1+i)^n )$ satisfy diagonal condition for some $n \geq 3$. Then, for all $\alpha \in \mathbb{Z}[i]$ s.t. $\gcd(\alpha, 1+i) = \delta$ where $\delta$ is unit, we have $\alpha^2 \equiv 1 \pmod{(1+i)^n}$. In $(\mathbb Z [i] /  (1+i)^3 )^\times$, we get $(-2+i)^2 = 4-4i-1 \equiv -1 \pmod{(1+i)^3}$ that makes $(-2+i)^4 \equiv 1 \pmod{(1+i)^3}$. Since $\mathbb Z [i] / ( (1+i)^n )$ satisfy diagonal condition for some $n \geq 3$, $(-2+i)^2 \equiv 1 \pmod{(1+i)^n}$ and hence $(-2+i)^2 \equiv 1 \pmod{(1+i)^3}$. Thus, $-1 \equiv 1 \pmod{(1+i)^3}$ so that $(1+i)^3 \mid -2$. This leads a contradiction, the statement we supposed is wrong. Hence, $\mathbb Z [i] / ( (1+i)^n )$ do not satisfy diagonal condition for all $n \geq 3$. 
	\end{proof}
	
	\vspace{2ex}\noindent{\textit{Remark 3.3} \ If the order of Gaussian integers modulo $\alpha$ is greater than or equal to $3$, then multiplication table of $\mathbb{Z}[i]/(\alpha)$ do not have diagonal condition.
	}
	
	\vspace{2ex}\noindent{\textbf{Lemma 3.4} \  $\mathbb Z [i] / (p)$ do not satisfy diagonal condition when $p$ are primes in $\mathbb Z$ s.t. $p \equiv 3 \pmod 4$.}
	\begin{proof} \ 
		Since $p$ are primes in $\mathbb Z$ such that $p \equiv 3 \pmod 4$, $p$ are Gaussian primes. Then $\mathbb Z [i] / (p)$ is a finite field which has order $p^2$ because $(p)$ is the maximal ideal. Thus, $(\mathbb Z [i] /  (p) )^\times$ is a cyclic group with its order is equal to $p^2 - 1 \geq 3$. By Remark 3.3, the statement of lemma is proved.
	\end{proof}
	
	\vspace{2ex}\noindent{\textbf{Theorem 3.5} \  $\mathbb Z [i] / ( (p)^n )$ do not satisfy diagonal condition when $n \geq 2$ and $p$ are primes in $\mathbb Z$ such that $p \equiv 3 \pmod 4$.}
	\begin{proof} \ 
		Assume that $\mathbb Z [i] / ( (p)^n )$ satisfy diagonal condition for some $n \geq 2$. Then, we have $\forall \alpha \in \mathbb Z [i]$ s.t. $\gcd(\alpha, p) = \delta$ where $\delta$ is unit, $\alpha^2 \equiv 1 \pmod{p^n}$. However, there exists $\beta \in \mathbb Z [i] / (p)$ such that the order of $\beta$ modulo $p$ is $p^2 - 1 \geq 3$. That makes $\beta^2 \not\equiv 1 \pmod p$, contradicts with $\beta^2 \equiv 1 \pmod p$ by assumption. Hence, the statement of theorem is confirmed.
	\end{proof}
	
	\vspace{2ex}\noindent{\textbf{Lemma 3.6} \ $\mathbb Z [i] / (\pi) $ and $\mathbb Z [i] / (\bar{\pi})$ do not satisfy diagonal condition when $\pi \bar{\pi} = p$ and $p$ are primes in $\mathbb Z$ such that $p \equiv 1 \pmod 4$.}
	\begin{proof} \ 
		By Theorem 2.3.7, $\pi \not\sim \bar{\pi}$ without loss of generality. Since $\pi$ is a Gaussian prime, $\mathbb Z [i] / (\pi) $ is a field because $(\pi)$ is a maximal ideal. Thus, $(\mathbb Z [i] / (\pi) )^\times$ is a cyclic group with order $p - 1$. Then, there exists $\alpha \in \mathbb Z [i]$ such that the order of $\alpha$ modulo $\pi$ is equal to $p - 1 \geq 3$ so that $\alpha^2 \not\equiv 1 \pmod \pi$. Hence, $\mathbb Z [i] / (\pi) $ does not satisfy diagonal condition.
	\end{proof}
	
	\vspace{2ex}\noindent{\textbf{Theorem 3.7} \  $\mathbb Z [i] / ( (\pi)^n )$ and $\mathbb Z [i] / ( (\bar{\pi})^n )$ do not satisfy diagonal condition when $n \geq 2$ and $\pi \bar{\pi} = p$ where $p$ are primes in $\mathbb Z$ such that $p \equiv 1 \pmod 4$.}
	\begin{proof} \
		By Theorem 2.3.7, $\pi \not\sim \bar{\pi}$ without loss of generality. Suppose $\mathbb Z [i] / ( (\pi)^n )$ satisfy diagonal condition. Then, $\alpha^2 \equiv 1 \pmod{(\pi)^n}$ for all $\alpha \in \mathbb Z [i]$. In $\mathbb Z [i] / (\pi)$, there exists $\beta \in \mathbb Z [i]$ such that the order of $\beta$ modulo $\pi$ is $p \geq 3$ so that $\beta^2 \not\equiv 1 \pmod \pi$. By the assumption, we need to have $\beta^2 \equiv 1 \pmod{(\pi)^n}$, a contradiction. Thus, the statement of theorem is proved.
	\end{proof}
	
	From Theorem 3.1 to Theorem 3.7 including remark and lemmas, we can generate the main theorem as followed:
	
	\vspace{2ex}\noindent{\textbf{Theorem 3.8} \ Multiplication table of $\mathbb Z [i] / (\alpha)$ have diagonal condition if and only if $\alpha = 1+i$ or $\alpha = (1+i)^2$.}
	
	\section{Open problems and interesting points}
	
	\qquad The author suggests readers to investigate diagonal condition on several interesting rings e.g. $\mathbb{Q}[i]$, \, $\mathbb{Z} \big[\frac{1+\sqrt{(-3)}}{2}\big]$, Real\slash \, Imaginary Quadratic fields ($\mathbb{Q}[\sqrt{D}]$ where $D$ is square-free, Real when $D>0$; Imaginary when $D<0$) or any other imaginary quadratic fields have units only $1$ and $-1$. \footnote[2]{This part is suggested by Prof. Jeffrey Clark Lagarias}Also, all related open problems such as Floor Quotient for the previous-mentioned domains with a variant problem analogs with what the author did in this paper.
	
		\begin{center}
			{\bf Disclosures}
		\end{center}
	\qquad The author declares that there are no conflicts of interest that could have influenced the objectivity of this research or the writing of this paper.
		
	
	\begin{center}
		{\bf Acknowledgements} 
	\end{center}
	
	\qquad The author would like to thank her BS advisor \textit{Assoc. Prof. Dr. Ajchara Harnchoowong} for inspiration on research of Number Theory. Additionally, the author would like to praise productive comments from \textit{Prof. Jeffrey Clark Lagarias} while attending the 2026 Gainesville international Number Theory conference at University of Florida, USA. The most important part to mention is her family, they always cheer her up in any stage of this research work.


\newpage
	
	\begin{thebibliography}{15}
		
		\bibitem[1]{Chebolu-24} Chebolu, S.K. ``What is special about the divisors of 24?." \textit{Math. Mag.} \textbf{85}, no. 5 (2012): 366--372.
		
		\bibitem[2]{Chebolu-12} Chebolu, S.K. and Mayers, M. ``What is special about the divisors of 12?." \textit{Math. Mag.} \textbf{86}, no. 2 (2013): 143--146.
		
		\bibitem[3]{K.Conrad} Conrad, K. ``The gaussian integers." Disponıvel em: \textit{http://www. math. uconn. edu/$\sim$ kconrad/blurbs/ugradnumthy/Zinotes.pdf}, Acesso em 28 (2016).
		
		\bibitem[4]{Dummit-Foote} Dummit, D.S. and Foote, R.M. \textit{Abstract algebra}. John Wiley \&
		Sons, Inc., Hoboken, NJ, third edition, 2004.
		
		\bibitem[5]{fromJeff-2} Edgar, G.A., Hensley, D. and West, D.B. ``Problems and Solutions." \textit{Amer. Math. Month.} \textbf{123}, no. 2 (2016): 197--204.
		
		\bibitem[6]{ContempAA} Gallian, J.A. \textit{Contemporary abstract algebra}. Chapman and Hall/CRC, 2021.
		
		\bibitem[7]{Sophie-Germain} Genzlinger, K. and Lockridge, K. ``Sophie Germain primes and involutions of $\mathbb{Z}_n^\times$." \textit{Inv., J. Math.} \textbf{8}, no. 4 (2015): 653--663.
		
		\bibitem[8]{fromJeff-1} Godinho, L., von Heymann, F. and Sabatini, S. ``12, 24 and beyond." \textit{Adv. Math.} \textbf{319} (2017): 472--521.
		
		\bibitem[9]{BookofThomas} Koshy, T. \textit{Elementary number theory with applications}. Academic press, 2002.
		
		
		\bibitem[10]{fromJeff-3} Poonen, B. and Rodriguez-Villegas, F. ``Lattice polygons and the number 12." \textit{Amer. Math. Month.} \textbf{107}, no. 3 (2000): 238--250.
		
	\end{thebibliography}
\end{document}